\documentclass[12pt]{scrartcl}
\usepackage{amssymb}
\usepackage{amsmath}
\usepackage{amsthm}
\usepackage{amssymb}
\usepackage{marginnote}
\usepackage{csquotes}

\begin{document}

\newcommand{\hs}{\hspace {8 cm}}
\newcommand{\cl}{\centerline}
\newcommand{\sms}{\smallskip}
\newcommand{\ms}{\medskip}
\newcommand{\bs}{\bigskip}
\newcommand{\R}{{\mathbf R}}
\newcommand{\C}{{\mathbf C}}

\def\al{{\alpha}}  \def\bet{{\beta}} \def\gam{{\gamma}}
 \def\del{{\delta}}
\def\eps{{\varepsilon}}
\def\kap{{\kappa}}                   \def\Chi{\text{X}}
\def\lam{{\lambda}}
 \def\sig{{\sigma}}  \def\vphi{{\varphi}} \def\om{{\omega}}
\def\Gam{{\Gamma}} \def\Om{{\Omega}}
\def\Del{{\Delta}}  \def\ups{{\upsilon}}
\def\th{{\theta}} \def\Th{{\Theta}}

\def\R{{\mathbf R}}  \def\Q{\mathbf {Q}} \def\P{{\mathbf P}}
\def\Z{{\mathbf Z}} \def\A{{\mathbf A}}
\def\N{{\mathbf N}}
\def\C{{\mathbf C}}  \def\V{{\mathbf V}}
\def\F{{\mathbf F}}
\def\W{{\mathbf W}}

\def\es{\emptyset}

\def\gom{{\mathfrak m}}
\def\gok{{\mathfrak k}}
\def\gog{{\mathfrak g}}
\def\goh{{\mathfrak h}}
\def\gou{{\mathfrak u}}
\def\gop{{\mathfrak p}}
\def\gob{{\mathfrak b}}
\def\goa{{\mathfrak a}}

\def\calA{{\mathcal A}}
\def\calP{{\mathcal P}}
\def\calO{{\mathcal O}}
\def\calB{{\mathcal B}}
\def\calT{{\mathcal T}}
\def\calM{{\mathcal M}}
\def\calS{{\mathcal S}}
\def\calF{{\mathcal F}}
\def\calX{{\mathcal X}}
\def\calY{{\mathcal Y}}
\def\calD{{\mathcal D}}
\def\calW{{\mathcal W}}
\def\calZ{{\mathcal Z}}
\def\calK{{\mathcal K}}
\def\calT{{\mathcal T}}
\def\calGr{{\mathcal Gr}}
\def\calHom{{\mathcal Hom}}

\def\Def{\textbf {Definition} }

\def\+{{\textbf{[P]}}}
\def\*{(*)}

\def\1{{\textbf 1.}}
\def\2{{\textbf 2.}}
\def\3{{\textbf 3.}}
\def\4{{\textbf 4.}}
\def\5{{\textbf 5.}}
\def\6{{\textbf 6.}}
\def\7{{\textbf 7.}}
\def\8{{\textbf 8.}}
\def\9{{\textbf 9.}}




\newcommand{\edoc}{\end{document}}

\def\al{{\alpha}}  \def\bet{{\beta}} \def\gam{{\gamma}}
 \def\del{{\delta}}
\def\eps{{\varepsilon}}
\def\kap{{\kappa}}                   \def\Chi{\text{X}}
\def\lam{{\lambda}}
 \def\sig{{\sigma}}  \def\vphi{{\varphi}} \def\om{{\omega}}
\def\Gam{{\Gamma}} \def\Om{{\Omega}}
\def\Del{{\Delta}}  \def\ups{{\upsilon}}
\def\th{{\theta}} \def\Th{{\Theta}}

\def\R{{\mathbf R}}  \def\Q{\mathbf {Q}} \def\P{{\mathbf P}}
\def\Z{{\mathbf Z}} \def\A{{\mathbf A}}
\def\N{{\mathbf N}}
\def\C{{\mathbf C}}  \def\V{{\mathbf V}}
\def\F{{\mathbf F}}
\def\W{{\mathbf W}}

\def\es{\emptyset}

\def\gom{{\mathfrak m}}
\def\gok{{\mathfrak k}}
\def\gog{{\mathfrak g}}
\def\goh{{\mathfrak h}}
\def\gou{{\mathfrak u}}
\def\gop{{\mathfrak p}}
\def\gob{{\mathfrak b}}
\def\goa{{\mathfrak a}}

\def\calA{{\mathcal A}}
\def\calC{{\mathcal C}}
\def\calP{{\mathcal P}}
\def\calO{{\mathcal O}}
\def\calB{{\mathcal B}}
\def\calT{{\mathcal T}}
\def\calM{{\mathcal M}}
\def\calS{{\mathcal S}}
\def\calF{{\mathcal F}}
\def\calX{{\mathcal X}}
\def\calY{{\mathcal Y}}
\def\calD{{\mathcal D}}
\def\calW{{\mathcal W}}
\def\calZ{{\mathcal Z}}
\def\calG{{\mathcal G}}

\def\calGr{{\mathcal Gr}}
\def\calHom{{\mathcal Hom}}
\def\itG{{\textit G}}
\def\Vect{{\mathrm{Vect}}}

\def\+{{\textbf[P]}}
\def\*{(*)}

\ms

\title{\large Stacks in Representation Theory \\ \vspace{5 mm}
{\small What is a continuous representation of an algebraic group ?}}

\author{Joseph Bernstein}

\maketitle

\ms

\section{{ Introduction.}}

In this note I would like to introduce a new approach to
(or rather a new language for) representation theory of groups.
Namely, I propose to consider a (complex) representation of a group $G$
as a sheaf on some geometric object.
    This point of view necessarily leads to a conclusion that
     the standard approach to (continuous)
representations of algebraic groups should be modified.

Let us start with a local or finite field $F$
and fix an algebraic group $\calG$ defined over $F$.
In the standard approach
we consider the set $G = \calG(F)$ of $F$-points of $\calG$
 as a topological group
and study an appropriate category $Rep(G)$ of continuous
representations of $G$.

\sms

The main goal of this note is to explain that
this approach is philosophically inconsistent.
In fact I will describe how to  extend the category $Rep(G)$ to
some larger category $\calM(\calG, F)$  that better
corresponds to our intuitive understanding of representations of $G$.

We will see that this category can be naturally described as a
product of categories $Rep(G_i)$ over all pure inner forms of the group $G$.
On the level of simple objects this means that
$Irr(\calM(\calG, F)) \approx \coprod Irr(G_i)$.
   This agrees with observation by several mathematicians
    (e.g.  by D. Vogan \cite{Vog})
   that  when we classify irreducible representations
   it is better to work with the union of sets $Irr(G_i)$ for several forms
   of the group $G$ than with one set $Irr(G)$.

\subsection{{ Representations and sheaves on Stacks. }}

   In order to describe the category $\calM(G)$  I propose to
    consider representations as sheaves on algebraic  stacks.

 Stacks play a more and more important role in
 contemporary Mathematics. Since they are not yet the common language in
 representation theory I will recall some basic notions related to  stacks.


\sms

   Informally, stack is a ``space" $X$ such that every point
   $x \in X$ is endowed with a group $G_x$ of automorphisms of
   inner degrees of freedom at this point.

   We see that in order to consider stacks  we should first fix a
   Geometric Environment,
   i.e. a category $\calS$ of spaces on which we
   model our stacks.
In fact the category $\calS$ should be considered together with some Grothendieck topology.
The standard term for such category $\calS$ is ``site".

\sms

    Usually one works with the following sites:

   \sms

   (i) Category of schemes over a field $F$ (with etale or smooth topology).

   (ii) Category of smooth manifolds (with usual topology).

   (iii) Category of locally compact Hausdorff topological spaces with usual 
   topology (or some natural subcategory of it, e.g. the category of totally disconnected spaces).

   (iv) Category $Sets$ of sets with discrete topology.

 \subsection{{  Groupoids. }}

       Stacks modeled on  the site $\calS = Sets$ are  {\textbf {groupoids}}. Let me remind that,
     by definition,  groupoid is       a category in which all morphisms are isomorphisms.

Groupoids represent rather elementary examples of stacks.  However they exhibit many
features of the general case. Also in this case the general ideas that I would like to explain
are much easier to understand.
For this reason I would like to discuss this case in some detail.

\subsubsection{{Examples of groupoids. }}

    To every discrete group  $G$  we assign
    the {\textbf{ basic groupoid}}  $BG = pt / G$ as follows:

\sms

   An object of the category $BG$ is a $G$-torsor $T$ (i.e. a non-empty $G$-set on which
$G$ acts transitively and free).
 The morphisms in this category are morphisms of $G$-sets.

\sms

   More generally, given an action of the group $G$ on a set $Z$ we define
   the {\textbf {quotient groupoid}} $BG(Z) = Z / G$ as follows:

   Object of $BG(Z)$ is a $G$-torsor $T$ equipped with a $G$-morphism $\nu: T \to Z$.
   Morphisms are morphisms of $G$-sets over $Z$.

  \subsubsection{   Representations as sheaves on groupoids.}

   Given a  groupoid $\calX$ it is natural to think about it
as a geometric object (some kind of a space). Then it is natural to consider sheaves on this space.

\sms

We define a {\textbf {sheaf}} $R$ (of complex vector spaces)  on  a groupoid $\calX$
 to be a functor $R: \calX \to Vect$, where $Vect$ is the category of complex
 vector spaces (we will see later why this notion is natural).
  We denote by $Sh(\calX)$ the category of sheaves on $\calX$.

  \sms

{\textbf {Claim.}}  The category $Sh(BG)$ is naturally equivalent to the category $Rep(G)$.

More generally, for every $G$-set $Z$ the category $Sh(BG(Z))$ is naturally
equivalent to the category $Sh_G(Z)$ of $G$-equivariant sheaves on $Z$
(see \ref{equiv-sheaves}).

\sms

This gives us a ``geometric" description of the category $Rep(G)$. This construction, that
is very elementary in case of groupoids, is the basis of the approach that I describe in this note.

\subsection{{ Topological groupoids.}}

 Let $G$ be a locally compact group. For technical reasons let us assume that it is
totally disconnected.  In this case we also can define the basic
groupoid $BG$ and quotient groupoids $BG(Z)$ as  stacks modeled on the site of
locally compact spaces. For any stack $\calX$
of this type we will define a category $Sh(\calX)$ of sheaves on $\calX$.
Using these constructions we can interpret the category
$Rep(G)$ as the category of sheaves on the stack $BG$.

\sms

  One of the ways to think about a stack modeled on the site of locally compact spaces
  is to interpret it as a topological groupoid.
 For example, using this interpretation
it is easy to show that for a $G$-space $Z$  the category $Sh(BG(Z))$ of sheaves
on the quotient stack
$BG(Z) = Z / G$ is naturally
equivalent to the category $Sh_G(Z)$ of $G$-equivariant sheaves
on $Z$ (see \ref{equiv-sheaves}).

     Technically,  working with topological groupoids is a little difficult.
      Later we  describe another way to define the category $Sh(\calX)$ of
sheaves on $\calX$  for stacks of this type that is technically simpler
(see section \ref{F-sheaves}).


\subsection{{ Algebraic groups and stacks.}}

Let us consider  more interesting case of an algebraic group $\calG$
 over a local (or finite) field  $F$.
 By analogy with the discrete case we define the {\textbf {basic
 stack}} $B\calG = pt/\calG$.
  This is an algebraic stack over the field $F$ (i.e. it is modeled on  the site $\calS$ of
  schemes over $F$ (see section \ref{stacks})).

\sms

  For any algebraic  stack $\calX$ over $F$ we will construct the category
 $Sh(\calX)$ of sheaves of complex vector spaces on $\calX$.
 The informal idea is that $F$-points
  $\calX(F)$ of the stack $\calX$ form a groupoid and  a sheaf $R$ on $\calX$
is just a sheaf $R$ on this groupoid.

   For finite fields this works fine. For local fields we should take into
 account the topology of the groupoid $\calX(F)$. \sms

\sms

  Given an algebraic group $\calG$ over $F$ we define
   the category $\calM = \calM(\calG, F)$  as the category $Sh(B\calG)$ of
  sheaves on the algebraic stack $B\calG$.
I call the objects of this category {\textbf {stacky $G$-modules}.}

\subsubsection{Two competing definitions}

Now starting with  algebraic group $\calG$ over $F$ we can consider
two competing definitions of a representation.

\sms

\Def \1 Category $Rep(G)$ \  obtained  from $\calG$ by a chain of constructions

\sms

\quad \quad  $\calG \quad \Longrightarrow$ \quad group  $ G = \calG(F) \quad \Longrightarrow$
\quad groupoid $\calY = BG$ $ \quad  \Longrightarrow \quad Sh(\calY)$

\ms

\Def \2 Category $\calM(\calG ,F)$ obtained from $\calG$ by a chain of constructions

\sms

\quad \quad   $\calG \quad \Longrightarrow$ \quad stack $\calX = B\calG
\quad  \Longrightarrow$ \quad groupoid
$\calX(F)$ \quad  $ \Longrightarrow \quad Sh(\calX(F))$

\sms

   The subtle point is that the groupoids $\calY$ and $\calX(F)$ are not
always equivalent.
So the category $\calM(\calG, F) = Sh(\calX(F))$ might be
not equivalent to the category $Rep(G) = Sh(\calY)$.

\sms

The standard notion of a continuous representation of $G$  is based on definition 1.
 \textbf{The main goal} of this note is to convince the reader that from 
 many points of view definition 2 is much more appropriate than the standard definition 1.

\subsection{{  Vogan's picture.}} \label{Vogan}

One advantage of this definition is that it gives an explanation to
representations that appear in Vogan's interpretation of the Langlands  correspondence.
Let me remind what is Vogan's suggestion (see \cite{Vog}, Conjecture 4.15 or \cite{ABV},
Theorem 1.18).

Vogan tried  to describe the Langlands correspondence in the following explicit way.
Let $F$ be a local field, $\calG$ a reductive algebraic group over $F$ and $G$
the group of its $F$-points. Consider on one side the set $Irr(G)$ of equivalence
classes of irreducible representations of $G$. On the other side consider the
set $Lan(\calG , F)$ of Langlands'
parameters $\Phi$ defined in terms of the dual group $L(\calG)$ and the field $F$.

Vogan tried to construct a canonical bijection between these two sets.
 He realized that in many
cases this can not work since the set $Irr(G)$ is just too small. But he also discovered that if
we replace this set by the disjoint union $\coprod Irr(G_i)$ of corresponding sets for all pure inner
forms $G_i$ of the group $G$ then this set has correct size (note that some of these forms
might be isomorphic -- then the ``same"  representation will appear in this list several times).
 In fact in many cases Vogan was
able to describe a bijection of this set with the set of Langlands' parameters.

\sms

In the language I propose  Vogan's conjecture can be formulated as follows.

 First of all, instead of the group $\calG$
we consider the algebraic stack $B\calG$. 
Instead of representations of the group $G = \calG(F)$ let us study the
category $\calM:= Sh(B\calG)$ of sheaves on this stack 
(i.e. the category of stacky $G$-modules).
Since for any pure inner form $\calG'$
of the group $\calG$ the stack $B\calG'$ is equivalent to $B\calG$,  the category
$\calM$ depends only on the pure inner class of $\calG$.

Our goal is to parameterize  the set $Irr(\calM)$ of isomorphism classes
of simple objects in $\calM$.
We will see in \ref{sheaves-algebraic-stacks} that  this set can be described
as a disjoint union over all pure inner forms $\calG_i$ of the group
$\calG$ of the sets $Irr(G_i)$, where $G_i = \calG_i(F)$.

Then the conjecture 4.15 in \cite{Vog} is essentially the statement
that the set  $Irr(\calM)$
 is in natural  bijection with the set $Lan(\calG, F)$ of Langlands' parameters.

 \sms

 {\textbf{Remarks}.}  \textbf{1.} In fact Vogan formulated his conjectures only for
 pure inner forms of quasi-split groups. Later they were generalized by Kaletha 
 to other forms (see \cite{Kal}).

 \textbf{2.} It would be interesting to understand whether one has
 a \textbf{canonical} bijection between these sets. 
 The constructions proposed in \cite{Vog}, \cite{ABV}
 are not quite canonical -- they depend on choice of some ``Whittaker data".

\smallskip

 I think that in fact they are not in canonical bijection. Namely,
 I think that the set of Langlands' parameters $Lan(\calG, F)$
 should be slightly modified -- and after this the resulting set $Lan'(\calG,F)$
 will be in a canonical bijection with the set $Irr(\calM)$.

 Possible suggestions for this modified set of Langlands' parameters
 can be found in \cite{BG} and sources listed there.

\subsection{About this note.}

In section \ref{sheaves} \ I remind basic facts about sheaves that are relevant
in Representation Theory.

In particular, I discuss the notion of an equivariant sheaf on a topological $G$-space $Z$.
This notion plays central role in any geometric approach to Representation Theory.
This is also an important computational tool -- in fact later on when I talk about
explicit description
of some object I mean a description in the language of equivariant sheaves.

In section \ref{example} \ I  describe a striking example that shows that
the standard definition of representations is not a good one.

   In subsequent sections I make various  comments on the notions we discussed above,
 indicate how to formulate precise technical definitions that describe these notions
 and describe some technical tools that help to make computations with them.

  In section \ref{groupoids} \ I  discuss the case of groupoids.
  In section \ref{stacks} \ I  shortly describe how to
  give a technical definition of a stack.

  In section \ref{F-sheaves} \  I
  discuss a technical definition of the category $Sh(\calX)$ of sheaves of
  complex vector spaces on
 a stack $\calX$ and describe how to make computations with them. In particular,
  I explain how one can describe the category $Sh(\calX)$ in
  terms of equivariant sheaves on topological spaces.

In section \ref{AG} \  I explain the algebro-geometric structure
that describes the relation between groupoids $\calX(F)$ for different fields $F$.

\ms

 This note is an expanded version of the lecture that I gave at
  the Fourth Conference of Tsinghua Sanya International Mathematics Forum (TSIMF)
   in December 2013. I would like to thank TSIMF organizers for the invitation.

\sms

  Main ideas about stacks in representations theory  came out of
  my research in the framework of the ERC grant 291612. My research was supported
  by the grant 533/14 of Israel Science Foundation.
  
  Much of my work on this subject was done during my visits to MPIM, Bonn.
  I would like to thank MPIM for very creative atmosphere.

  I would like to thank A. Vistoli for his remarks about stacks that helped me to
understand this notion. I also thank D. Zagier for critical remarks about an
earlier version of this note.

  I would like to thank the referee for several useful suggestions.


\section{ Sheaves and equivariant sheaves.}   \label{sheaves}

\subsection{Sheaves relevant in representation theory.}

Let $F$ be a finite or local field, $\calG$ an algebraic group over $F$.
We consider the group $G = \calG(F)$ of its $F$-points as a topological group.

In order to study representations of the group $G$ we usually place this group in some
geometric environment.
Namely we fix a site $\calS$ appropriate for the group $G$ and for every space $Z \in \calS$
we consider some category $Sh(Z)$ of sheaves of complex vector spaces  on $Z$.
 \sms
  Let us describe this in more detail.

\sms

(i)  Let $F$ be a finite field. We consider the site $\calS$ of finite discrete sets.
 For a space $Z \in \calS$ we denote by $Sh(Z)$  the category of all sheaves of
 complex vector spaces on $Z$.

\sms

(ii) Let $F$ be a local non-Archimedean field. We consider the site $\calS$ that consists
of Hausdorff locally compact totally disconnected spaces $Z$ with countable base of open subsets
($l$-spaces in terminology of \cite{BZ}).

Given an $l$-space $Z$ we consider the category $Sh(Z)$ of all sheaves of complex vector
spaces on $Z$ (in \cite{BZ} they were called $l$-sheaves).

\sms

(iii) Let $F$ be a local Archimedean field, i.e. $F = \R$ or $F = \C$. In this case I see several
candidates
for the site $\calS$.

    We can consider $\calS$ to be the category of smooth manifolds.
For a manifold $Z$ the category $Sh(Z)$ is the category of sheaves of $\calO_Z$-modules,
where $\calO_Z$ is the sheaf of smooth complex valued functions on $Z$.

\sms

Another possibility is to work with the site of Nash manifolds and some sheaves of Schwartz functions
(or distributions) on these manifolds (see \cite{AG}).

I do not know what is the correct approach to this case. So in what
 follows I mostly deal  with finite and local non-Archimedean fields.

\ms

{\textbf {Remarks}.} \1 In this note for simplicity I discuss only representations over the field $k = \C$
of complex numbers. In case when $F$ is a finite or non-Archimedean field we can consider other
fields  $k$ of coefficients.

If $F$ is a local non-Archimedean field and $k$ is an extension of $F$ we can try to include the
theory of  locally analytic representations by considering the site $\calS$ of analytic
 manifolds over $F$
and categories $Sh(Z)$ of sheaves of modules over the sheaf $\calO_Z$
of locally analytic functions on $Z$.
It seems that the stacky language that I will introduce might be applicable
and  useful in this theory.

\sms

\2  The structures I described are important when we are trying to specify the category
$Rep(G)$ of representations of the group $G$ that we would like to study.

In case when the group $G$ is defined over a finite or local non-Archimedean field $F$
there is a consensus what is the ``correct" category of representations $Rep(G)$.
 Namely, in case of a finite field we consider the category of all representations and
in case of a local non-Archimedean field we consider the category of smooth
 representations (smooth means that every vector has open stabilizer in the group $G$).

  The situation for  real groups is different -- the correct choice of the appropriate
  category $Rep(G)$ is a very
non-trivial question. In case of reductive groups this was done by Casselman-Wallach
(see \cite{Cass},\cite{W}, \cite{BK}).
For general real algebraic group I do not know a good candidate for this category.

 In this note I am not going to discuss this tricky question. For this reason I will mostly deal with finite
 and local non-Archimedean fields.

\subsection{Equivariant sheaves.} \label{equiv-sheaves}

From now on we assume that all topological spaces we consider are $l$-spaces
(locally compact, totally disconnected with countable base of open subsets).

   Let $G$ be a topological group and $Z$ a $G$-space, i.e. $Z$ is equipped with
a continuous action $a: G \times Z \to Z$.
We denote by $Sh_G(Z)$ the category
 of $G$-equivariant sheaves (of complex vector spaces)  on $Z$.
This category will play central role in what follows.

\sms

Recall, that a $G$-equivariant sheaf is a sheaf $R$ on $Z$ equipped with
an isomorphism $\al: a^*(R) \to pr_Z^*(R)$ of two liftings of $R$ to the
space $G \times Z$ satisfying some natural cocycle condition (this condition is that after the lifting
to the space $G \times G \times Z$ two morphisms  $\nu, \mu: (m \cdot a)^*(R) \to pr^*(R)$ of
sheaves naturally constructed from $\al$
should coincide (see details of the  definition and discussion in \cite{BL}, Part1, section $0$).

\sms

{\textbf {Remark}.} In case of real groups we can work with the site $\calS$ of
smooth manifolds and with sheaves of $\calO$-modules. In this case
  the definition of equivariant sheaves formally  looks exactly the same, but has
quite different  geometric meaning.
The reason is that the pullback functor $a^*$ in the category of $\calO$-modules
 is quite different from the pullback functor in the category of sheaves.

\ms

   Let me remind  two  standard facts about equivariant sheaves
 (here we assume that $G$ is an $l$-group).

\sms

 {\textbf {Fact 1}.} Let $Z$ be a point. Then the category $Sh_G(Z)$ is equivalent
 to the category $Rep(G)$ of {\textbf {smooth}} representations of $G$.

 \sms

{\textbf {Fact 2}.} Suppose that $Z$ is a quotient space of the group $G$.
 Fix a point $z \in Z$ and denote by $H$ its stabilizer in $G$. Then we have natural
equivalences of categories
     $Sh_G(Z) \approx Sh_H(z) \approx Rep(H)$.

\section{An example.}  \label{example}

In this section I describe a striking example that illustrates what is wrong with
the standard approach.

\subsection{ Heuristics.} My example  is based on the following  heuristic geometric principle.

Let $a: \calG \times \calZ \to \calZ$ be a transitive action of an algebraic
group $\calG$ on an algebraic variety $\calZ$.
Passing to $F$-points we get a continuous action $a: G \times Z \to Z$.

This action is usually not transitive. In many cases we can write $Z$ as a
union of open orbits $Z = \bigsqcup Z_i$, $i = 1,...,n$.

\sms

{\textbf {Heuristic Geometric Principle}.}

\sms

1. The space $Z$ is ``good", i.e. it is easy to describe.

\sms

2. Every individual orbit $Z_i$ might be  a ``bad" space, that
means that it is difficult to describe.

\sms

{\textbf {Illustration}.} Consider the space $V$ of real symmetric  $8 \times 8$ matrices and
try to give explicit descriptions of subsets $Z, Z_3  \subset V$ that describe non-degenerate
quadratic forms and quadratic forms of signature $3$ on $\R^8$ respectively.

\subsection{{ An example -- representations of orthogonal groups.}}

Fix an  $n$-dimensional vector space $V$ over $F$.

 The group $G = GL(V )$ acts on the space $Z$ of non-degenerate quadratic forms.

Fix a form $Q \in Z$ and denote by  $H$ its stabilizer in $G$ (i.e. $H$ is the  orthogonal group $O(Q)$).
Let us denote  by $Z_Q$ the $G$-orbit of $Q$ in $Z$. Then we have an equivalence of categories
$Rep(H) \approx Sh_G(Z_Q)$.

\ms

    According to the heuristic geometric principle the space $Z_Q$ might be (and often is)
a ``bad" space. This means that the  category $Rep(H) \approx Sh_G(Z_Q)$ might be a
 ``bad" category (i.e. very difficult to describe).
  In other words, the ``natural" problem of classification of irreducible representations
of the orthogonal group $H = O(Q)$ turns out to be not that natural.

 \sms
     However we see that the bad category $Rep(H) \approx Sh_G(Z_Q)$  can be naturally
extended to a larger  category
     $\calM := Sh_G(Z)$ of all  $G$-equivariant sheaves on the good algebraic space $Z$.
 We can expect (and this is really the case) that this larger category $\calM$ is a ``good" category.

\sms

     Note that we have a natural decomposition $\calM \approx \prod  Sh_G(Z_i)$,
where $Z_i$ are $G$-orbits  in $Z$. In particular the set $Irr(\calM)$ of isomorphism classes of simple
objects in $\calM$ is a disjoint union of sets $Irr(Sh_G(Z_i))$.

     It seems reasonable to assume that the classification of
     simple objects of the category
     $\calM$ might be relatively simple  problem, but then to sort out
     which of them are related
     to the orbit $Z_Q$ might turn out to be much more difficult problem
     (and it is not clear whether this problem is a meaningful one).

\ms

 The example we are considering suggests a certain pattern  that seems to work
also  in the general case.
 Namely, we see that if $G$ is a group of points of an algebraic group
then the category $Rep(G)$ of its representations might be a bad category,
but we can include it as a direct factor  into some larger good
category $\calM(G)$.

\sms

I had this example in mind for some time until I realized how one can define
this larger category using sheaves on stacks.

\section{{ Some remarks about groupoids}} \label{groupoids}

\subsection{{Equivalence of groupoids}}

We know that if two objects of some category are isomorphic
then it is better to consider them as two realizations of the same geometric structure.
Similarly,
if two groupoids $\calX$ and $ \calY$ are equivalent (as categories) we
 can assume that they represent
two realizations of the same geometric structure.
 \sms
\sms

A subtle point here is that the equivalences between these groupoids form a groupoid.
This means that if we fix an equivalence $Q: \calX \to \calY$ then this equivalence itself has
automorphisms, and it is not immediately clear how we should think about them.

\sms

  Similarly, if we would like to show that two groupoids are {\textbf {canonically}}
equivalent we have to construct an equivalence between them and show that this equivalence
 is defined up to a canonical isomorphism.

\sms

{\textbf {Example}.} Consider an action $a: G \times Z \to Z$.

Let us define a groupoid $BG_0(Z)$ as follows:

\sms

  $ Ob(BG_0(Z)) = Z, Mor(BG_0(Z)) = G \times Z$, where morphism $(g,z)$ is a morphism from the object $z$ to the object $gz$.

\sms

The groupoid $BG_0(pt)$ we denote by $BG_0$. \sms

\sms

{\textbf {Claim}.} The groupoid $BG_0(Z)$ is canonically equivalent to the groupoid $BG(Z)$.

\sms

The groupoid $BG_0(Z)$ might be considered as a
``matrix" version of the groupoid $BG(Z)$. It is better suited for computations.

\sms


\subsection{{  Theory of  groups and theory of groupoids.}}

I would like to explain that the theories describing groups and groupoids
are essentially equivalent.

\ms

{\textbf {Proposition}.} \1 Every groupoid $\calX$ is canonically decomposed
as a disjoint union of connected groupoids
(groupoid is connected if all its objects are isomorphic).

\2 A connected groupoid $\calY$ is equivalent to the basic
groupoid for some group $G$.

\sms

 This result shows that any question about groupoids can be reduced
 to a question in group theory.
In other words, the difference between theories of groups and groupoids
  is in their emphasis.

  In my opinion the relation between the theory of groupoids and the group theory
is very similar to the relation between  linear algebra and matrix calculus.
While these two theories are basically equivalent, clearly linear algebra is much more intuitive.
So I expect that eventually the stacky approach will become a standard tool in representation theory.

\sms

\subsubsection{{   Equivalence between groups and connected groupoids.}}

 The group $G$ corresponding to a connected groupoid $\calY$
  is not defined canonically.
   It depends on a choice of an object $Y \in \calY$.

   Namely, given an object $Y$ we can define a group $G = G_Y$ by $G:= Aut(Y)$. 
   Then we get canonical equivalence of categories $Q = Q_Y: \calY \to BG$, 
   defined by $Q(X) = Mor(Y,X)$.
   
   If we pick another object $Y'$ we get a different group $G'$ and a 
   different equivalence $Q': \calY \to BG'$.

   Note that any choice of an isomorphism  $\nu: Y \to Y'$
 defines natural  isomorphisms $G \backsimeq G'$
   and  $Q \backsimeq Q'$. However there is no preferred  choice for such an isomorphism $\nu$.

   \subsection{Examples of groupoids.}

   The next three constructions show that  in Mathematics we usually
    encounter groupoids and not groups.

\subsubsection{{  Multiplicative groupoid of a category.}}

{\textbf {Construction I.}} \  Starting with any category $C$ we construct the
 {\textbf {multiplicative groupoid}} $C^* = Iso(C)$
 that has the same collection of
objects as category $C$ and isomorphisms of $C$ as morphisms.
 \sms

\ms

{\textbf {Example 1.}} $C = Finsets$ --  the category of finite sets.

\sms

In this case the groupoid $Iso(C)$ is essentially the collection of all
symmetric groups $S_n$.

 \ms

 {\textbf {Example 2.}} $C = Vect_k$ -- the category of finite dimensional vector
spaces over a field $k$.

\sms

In this case the groupoid  $Iso(C)$ describes the collection of
groups $GL(n,k)$ for all $n$.

 \subsubsection{{  Poincar\'e groupoid.}}

  {\textbf {Construction II.}} \  Poincar\'e groupoid  $Poin(X)$ of a topological space $X$.

 \sms

 Objects of $Poin(X)$ are points of $X$. Morphisms from $x$ to $y$ are homotopy
classes of paths from $x$ to $y$.

 \sms

 If the space $X$ is path connected then the groupoid $Poin(X)$ is connected.
 For any point $x \in X$ the group
 $Aut_{Poin(X)}(x)$ is the fundamental group $\pi_1(X,x)$.

 This shows that the Poincar\'e groupoid is more basic notion than the fundamental group.

\sms

\subsubsection{{  Galois groupoid.}}

{\textbf  {Construction III. }} \  Galois groupoid $Gal(F)$  of a field $F$.

 \sms

 Objects of the groupoid $Gal(F)$ are field extensions $F \to \Om$ such that
 $\Om$ is an algebraic closure of $F$.
 Morphisms are isomorphisms of field extensions.

 \sms

   The groupoid $Gal(F)$  is connected.
 If we fix an algebraic closure $\Om$ then by definition the group $Aut_{Gal(F)}(\Om)$
is the absolute  Galois group  $Gal(\Om / F)$.

   Again  we see that the notion of Galois groupoid is more basic than
the notion of Galois group.

  \sms

   Note that the constructions of Poincar\'e and Galois groupoids are very similar.

\sms

\section{  What is a stack ?} \label{stacks}

   Let us fix some site $\calS$. I would like to describe
   the notion of a stack $\calX$ modeled on $\calS$.
   I assume two features of this notion.

\sms
 \sms
1. For every two stacks $\calX, \calY$ the collection of morphisms
 from $\calX$ to $\calY$
forms a groupoid  $Mor(\calX,\calY)$.

2. Every object $S \in \calS$ is a stack.

\sms

      The natural idea is to characterize a stack $\calX$ by the collection of groupoids
      $\calX(S):= Mor(S, \calX)$ for all objects $S \in \calS$.

      \sms

      In fact usually it is enough to know the groupoids $\calX(S)$
      for objects $S$ in some  subcategory  $\calB \subset \calS$ provided it is large enough.
      For example, if $\calS$ is the category of schemes we can restrict everything to the
      subcategory $\calB$ of affine schemes.

\subsection{{ Informal technical definition of a stack.}}

 Fix a large subcategory $\calB \subset \calS$. We define a stack $\calX$ over the
 site $\calS$  to be the following collection of data:

\sms

(i) To every object $S \in \calB$ we assign a groupoid $\calX(S)$

(ii) To every morphism  $\nu: S \to S'$ in $\calB$  we assign a functor
$\calX(S') \to \calX(S)$

(iii) To every composition of morphisms in $\calB$ we assign an
isomorphism of the corresponding functors.

\sms

 This data should satisfy a variety of compatibility conditions . These include

 (i) Compatibility conditions for isomorphisms we have chosen.

  (ii) Descent properties for
morphisms and for objects with respect to the Grothendieck topology on the site $\calS$.

(iii) Some finiteness conditions.

(iv) We also usually assume that the stack $\calX$ is
dominated by some object $Z \in \calS$ (for example the quotient stack $Z / G$ that
we will describe bellow is dominated by an object $Z$).

\sms

 A relatively elementary exposition of stacks one can find in \cite{Fan}.
 More detailed and more sophisticated exposition see in \cite{Vis}.

\subsection{{Stacks modeled on the site $Sets$}}

Consider the site $\calS = Sets$. To describe a stack $\calY$ modeled on $\calS$ we have to
assign to every set $S \in \calS$ a groupoid $\calY(S)$.

\sms

{\textbf{Example.}} Given a group $G \in \calS$ we construct the basic stack $BG$ as follows:

For any set $S$ the objects of the groupoid $BG(S)$ are principle $G$-bundles $P$ over $S$
and morphisms are $G$-morphisms over $S$.

 Let me remind that principle $G$-bundle over $S$
is a $G$-set $P$ equipped with a $G$-morphism $p: P \to S$, where $G$ acts trivially on $S$,
that locally on $S$ is isomorphic to a trivial $G$-bundle $pr: G \times U \to U$.

\sms

More generally, if $Z$ is a $G$-space we define the quotient stack $BG(Z)$ as follows:

For every set $S$ the objects of the groupoid $BG(Z)(S)$ are principle $G$-bundles over $S$ equipped with
a $G$-morphism to $Z$ and morphisms are  morphisms of $G$-bundles over $Z$.

\ms

Note that in case of sets we can restrict everything to a subcategory $B \subset \calS$ that
contains just one object $pt$ -- this is a big enough subcategory. This shows that
every stack $\calY$ modeled
on the site $\calS$ can be completely described by the groupoid $\calX = \calY(pt)$.
This explains why the stacks in this case can be considered as groupoids.

\subsection{Stacks modeled on the site $S$ of $l$-spaces.}

Consider the site $\calS$ of $l$-spaces. Let $G \in \calS$ be an $l$-group.
We define the basic stack $BG$ as follows:

For every $S \in \calS$ the objects of the groupoid $BG(S)$ are principal $G$-bundles over $S$ and
morphisms are morphisms of $G$-bundles.
In this case again a $G$-bundle is  a morphism of $G$-spaces $p:P \to S$ that locally in $S$ is
isomorphic to the trivial $G$-bundle $pr: G \times U \to U$.

\sms

 Similarly one defines a stack $BG(Z) = Z/G$
for a $G$-space $Z$.

\sms

In case of the site of $l$-spaces we can restrict everything to the
subcategory $B \subset \calS$ consisting
of compact spaces. Using this fact it is easy to check that the basic stack $BG$
can be described in terms of a topological groupoid $\calX = BG(pt)$.

Up to equivalence this groupoid can be explicitly described as follows: \
 $\calX$ has one object $X$  and its automorphism group $Aut(X)$ is
the topological group $G$.

\subsection{Stacks modeled on the site of $F$-schemes.}

\subsubsection{Basic stack $B\calG$.}

Let us fix a field $F$ and consider the site $\calS$ of schemes over $F$.
 Let $\calG \in \calS$ be an algebraic group defined over $F$.
We define  the basic stack $B\calG$
 in the same way as before.

Namely, for an affine $F$-scheme $S$ an object of the  groupoid $B\calG(S)$
is a principal $G$-bundle $P$ over $S$. Morphisms are morphisms of $G$-bundles.

\sms

In this case the notion of a principal $G$-bundle is more subtle. Namely, we should
 consider morphisms $p: P \to S$ of $G$-schemes that should be locally trivial in etale topology.
This is more sophisticated notion.

\subsubsection{Torsors.}

Let us consider in more detail the case of a point, i.e. assume that $S = Spec(F)$.
In this case we have to consider a morphism $p: T \to Spec(F)$, where $T$ is an $F$-scheme
 with an action of the group $G$
 that is locally trivial in etale topology (such $T$ is called a \textbf {$G$-torsor}).

 \sms

 The condition of local triviality is equivalent to the statement that after the base change
 from $S = Spec(F)$ to $S' = Spec(F')$, where $F'$ is an algebraic closure of $F$, the torsor $T'$
 will be isomorphic as a $G'$-space to the trivial $G'$ torsor $G'$.

 \sms

 Here is a typical example of torsors.

 \sms

 {\textbf{Example.}} Let $K$ be a finite separable field extension of $F$.
  Multiplicative groups $K^*$ and $F^*$ we consider as $F$-points of 
   algebraic groups $\calK^*, \calF^*$ defined over $F$. Then we have the norm 
   morphism of algebraic groups
 $N: \calK^* \to \calF^*$. We denote
   by $G$ the kernel of this morphism.
  This is an algebraic group defined over $F$.

  For any point $a \in F^*$ we consider the variety $\calT_a = N^{-1}(a) \subset \calK^*$.
  This variety is defined over $F$ and  is a $G$-torsor.

  Two torsors $\calT_a, \calT_b$ are isomorphic iff the quotient  $a/b$ lies in the image of  $K^*$ under the norm map
  $N: K^* \to F^*$.

\subsubsection{Quotient stack $B\calG(\calZ)$.}

\sms

    Let  $\calG$ act on a scheme $\calZ \in \calS$. \
 We define the quotient stack
    $\calX = B\calG(\calZ) = \calZ / \calG$ as follows:

    \sms

    For an affine scheme $S \in \calS$ the category $B\calG(Z)(S)$ consists of principal $G$-bundles
    $p: P \to S$ equipped with a $G$-morphism $ P \to Z$ and morphisms are morphisms of
    principal $G$-bundles over $Z$.

\sms

\section{{  Sheaves on stacks.  }} \label{F-sheaves}

\subsection{{  Technical definition of sheaves on stacks.}}

Let us fix a site $\calS$. In situations we consider we have the notion of sheaves
(of complex vector spaces) on spaces $Z \in \calS$. In other word, to every space
$Z \in \calS$ we assign a category $Sh(Z)$, every morphism $\nu: Z \to W$ in
$\calS$ induces a functor $\nu^*: Sh(W) \to Sh(Z)$ with a compatibility isomorphisms
for the products of two morphisms.

\ms

   {\textbf{Basic examples.}} \1 Let $S$ be the site of $l$-spaces. To every space $Z \in \calS$
we assign the category $Sh(Z)$ of $l$-sheaves (i.e sheaves  of complex vector spaces) on $Z$;
the functor  $\nu^*$ is the usual pullback functor.

\sms

\2 Fix a local non-Archimedean field $F$ and consider the site $\calS$ of $F$-schemes of finite type.
For every scheme $Z \in \calS$ we consider the set $Z(F)$ of its $F$-points as an $l$-space
and we set $Sh(Z): = Sh(Z(F))$ (i.e. a sheaf $R$ on $Z$ is just a sheaf of complex vector spaces
on the $l$-space $Z(F)$).

    \sms

We would like to extend these categories of sheaves  to stacks modeled on the site $\calS$.
In other words, given a stack $\calX$ we would like in some natural way
to assign to it a category $Sh(\calX)$   so that for a space $Z \in  \calS$
considered as a stack it will be the category $Sh(Z)$. Let me indicate some general
strategy how to construct the category $Sh(\calX)$.

\ms

 Suppose we already have some notion of sheaves on stacks.
    Fix a sheaf $R$ on a stack $\calX$.
    Then for any  space $S \in \calS$ and any point $p \in \calX(S) = Mor(S, \calX)$
    we get a sheaf $R_p = p^*(F) \in Sh(S)$.
     We also get a family of isomorphisms connecting these sheaves.
    Now we can try to use these sheaves and isomorphisms to
   characterize the sheaf $R$ as follows.

\sms

    {\textbf {Informal definition}.}   A sheaf $R$ on the stack $\calX$ is a collection of
 sheaves $R_p \in Sh(S)$
   for all spaces $S \in \calS$ and all
      morphisms $p: S \to \calX$ and a collection of isomorphisms satisfying correct
      compatibility relations.

\sms

{\textbf{Claim.}}\1  Let $\calS$ be the site of sets, $\calY$ a stack modeled on $\calS$.
Then the category $Sh(\calY)$ is naturally equivalent to the category $Sh(\calX) := Funct(\calX, Vect)$,
where $\calX$ is the groupoid $\calY(pt)$ that characterizes stack $ \calY$.

\sms

\2 Let $\calS$ be the site of $l$-spaces. Consider an action of an $l$-group $G$ on an
$l$-space $Z$ and denote by $BG(Z)$ the quotient stack, Then the category $Sh(BG(Z))$ is
naturally equivalent to the category $Sh_G(Z)$ of $G$-equivariant sheaves on $Z$.

\ms

{\textbf{Remark}.} In some cases the pullback functor is defined only partially (for some good morphisms).
This happens, for example, if we consider the site of manifolds and sheaves of distributions on them.
In this case pullback functor is well defined only for submersions.

In these cases we often can define appropriate notion of a sheaf for a stack $\calX$ that is dominated by an
object $Z \in \calS$, i.e. has a morphism $p: Z \to \calX$, provided that this morphism $p$ is good.
The quotient stacks that we are interested in usually have this property.

\subsection{{ How to describe sheaves on an algebraic stack. }} \label{sheaves-algebraic-stacks}

    Let $\calX$ be an algebraic stack over $F$. I would like to give a
    convenient description of
    sheaves on the stack $\calX$ in terms of equivariant sheaves.
I will do this for the case of a quotient stack $\calX \approx \calZ/\calG$.

    \sms

    Let $T_1,...T_r$ be representatives of isomorphisms classes of $\calG$-torsors.
     They are described by elements of $H^1(Gal(F), \calG(\bar{F}))$ (for simplicity we assume that this
set is finite; this is always the case when $char(F) = 0$).

    For every index $i$ consider the group  $\calG_i = Aut_\calG(T_i)$ -- this gives us
 the collection of all   pure inner forms of the group $G$. Also we consider a $\calG_i$-scheme
 $\calZ_i = Mor_\calG(T_i, \calZ)$.

 Passing to $F$ points we construct a topological group $G_i = \calG_i(F)$
 and a topological $G_i$-space $Z_i = \calZ_i(F)$.

  \sms

    {\textbf {Claim}.} The category $Sh(\calX)$ of sheaves on the stack $\calX$
     is naturally equivalent to $ \prod Sh_{G_i}(Z_i)$ (product of categories of equivariant sheaves).

\sms

    In particular we see that the collection of simple
    objects of the category $Sh(\calX)$ is a disjoint union of
    collections of simple $G_i$-equivariant sheaves on $Z_i$.

\sms

{\textbf{Remark}.} In case when $\calZ = pt$ we see that $Irr(\calM(\calG,F)) = \coprod Irr(G_i)$.
   This means that the set $Irr(\calM(\calG, F))$ of equivalence classes of irreducible stacky $G$-modules
   can be naturally described as the disjoint union  $ \coprod Irr(G_i)$ taken over all pure inner forms
   $\calG_i$ of the group $\calG$.

\subsection{{  Reduction to the case of the group $GL(n)$.}} \label{GL(n)}

Let me present one more  description of the category of
 sheaves on a quotient stack $\calX = \calZ / \calG$
 that is often convenient in computations.

\sms

{\textbf {Construction}}. \  Suppose that $\calG$ is a linear algebraic group.
Then we can imbed it into a group
$\calP$ isomorphic to $GL(n)$.

\sms

Using this we can realize our quotient stack $\calX = \calZ / \calG$
as the quotient stack  $\calW / \calP$, where  $\calW = \calP \times_{\calG} \calZ$.

\ms

The group $\calP$ has only one pure inner form
 (this is Hilbert 90 theorem). This implies that
 the category  $Sh(\calX)$ can be realized as the category
 $Sh_P(W)$ of $P$-equivariant sheaves on $W$, where $P = \calP(F)$
 and $W = \calW(F)$.

\sms

{\textbf{Remark}.} The example in section \ref{example} is a special case
of this construction.

\sms

\section{{  Algebro-geometric structure of the groupoids $\calX(F)$.}} \label{AG}

   Let $\calX$ be an algebraic stack.
Am important role in the study of sheaves on the stack
$\calX$ plays the fact that the collection of groupoids
$\calX(F)$ for different fields $F$ has an algebro-geometric structure.

\ms

{\textbf {Proposition}.}   Let $L \supset F$ be a finite Galois field extension and
$\Gamma = Gal(L / F)$ its Galois group. Then the group $\Gam$ acts on the
groupoid $\calX(L)$ and the {\textbf {fixed point groupoid}} $\calX(L)^\Gamma$ is
naturally equivalent to the groupoid $\calX(F)$.

\ms

Here the fixed point groupoid $\calX(L)^\Gam$ is defined in a standard categorical
manner. Its object is an object  $X$ of the groupoid $\calX(L)$ equipped with
a collection of isomorphisms $\al_\gam: X  \to \gam(X)$
satisfying natural compatibility conditions.

\ms

This proposition is just a reformulation of the descent property for the stack $\calX$.



\begin{thebibliography}{BBB}

\bibitem[ABV]{ABV} J. Adams, D.Barbash, D. Vogan, The Langlands Classification and Irreducible Characters
for Real Reguctive Groups, Progress in Mathematics v. 104(1992), Springer

\bibitem[AG]{AG} A. Aizenbud, D. Gourevitch, Schwartz functions on Nash manifolds,
 International Math. Research Notices (2008), no.5

\bibitem[BK]{BK} J. Bernstein, B. Kroetz, Smooth Frech\'et representations of Harish-Chandra modules,
Israel J. Math. 199 (2014), no.1,  45-111

\bibitem[BL]{BL} J. Bernstein, V. Lunts, Equivariant sheaves and functors,
Lecture Notes in Mathematics, v.1578 (1994).


\bibitem[BZ]{BZ} J. Bernstein, A. Zelevinski, Representations of the group
$GL(n,F)$, where $F$ is a
non-Archimedean local field, Russian Math. Surveys, 31-3 (1976),  1-68


\bibitem[BG]{BG} K.Buzzard, T. Gee, The conjectural connections between Automorphic Representations and Galois Representations, Proceedings of the ELM Durham Symposium, 2011 and arxiv:1009.0785.

\bibitem[Cass]{Cass} W. Casselman, Canonical extensions of Harish-Chandra modules to
 Representations of $G$,
Canad. J. Math. 41(1989), no.3, 385-438


\bibitem[Fan]{Fan} B. Fantechi, Stacks for Everybody,
European Congress of Mathematics,
Progress in Mathematics, Volume 201 (2001),  349-359



\bibitem[Kal]{Kal}  T. Kaletha, Rigid Inner Forms of Real and $p$-adid Groups,
 arxiv: 1304.3292

\bibitem[Vis]{Vis}  A. Vistoli,  Grothendieck topologies,
fibered categories and descent theory, Fundamental Algebraic Geometry,
1-104, Math. Survey Monogr.  123 (2005), AMS.


\bibitem[Vog]{Vog}  D. Vogan,  The local Langlands conjecture, Representation theory of groups and algebras,
Contemporary Mathematics, 145(1993), AMS.

\bibitem[W]{W} N. Wallach,, Real reductive groups II (section 11), Academic Press, 1992









\end{thebibliography}
\end{document}